\documentclass[reqno, 10pt]{amsart}

\usepackage{amsfonts}
\usepackage{graphicx}

\newtheorem{theorem}{Theorem}[section]                                          
\newtheorem{proposition}[theorem]{Proposition}                          
\newtheorem{lemma}[theorem]{Lemma}

\newcommand{\no}{\noindent}
\newcommand{\nn}{\nonumber}
\newcommand{\oo}{\infty}

\newcommand{\ra}{\rightarrow}

\newcommand{\<}{\langle}
\newcommand{\N}{\mathbb{N}}
\newcommand{\Z}{\mathbb{Z}}

\makeatletter
\@addtoreset{equation}{section}
\makeatother

\bibdata{prob}
\bibstyle{alpha} 

\title{Hydrodynamics for a non-conservative Interacting Particle System}

\author{Glauco Valle}

\date{}

\address{
\newline
UFRJ - Departamento de m\'etodos estat\'{\i}sticos do Instituto de Matem\'atica.
\newline  Caixa Postal 68530, 21945-970, Rio de Janeiro, Brasil
\newline
e-mail:  \rm \texttt{glauco.valle@dme.ufrj.br}
}

\thanks{Work supported by CNPq and FAPERJ.}

\begin{document}
  
\maketitle

\begin{abstract}
We obtain the hydrodynamic limit of one-dimensional interacting particle systems describing the macroscopic evolution of the density of mass in infinite volume from the microscopic dynamics. The processes are weak pertubations of the symmetric exclusion process by shift operators describing the spread of particles around positions where new sites are created.
\end{abstract}

%%%%%%%%%%%%%%%%%%%%%%%%%%%%%%%%%%%%%%%%%%%%%%%%%%%%%%%%%%%%%%%
\section{Introduction} 
%%%%%%%%%%%%%%%%%%%%%%%%%%%%%%%%%%%%%%%%%%%%%%%%%%%%%%%%%%%%%%%
\label{sec:intro}

The basic model of interacting particle systems that gives a rough microscopic description of the evolution of the mass density profile of an incompressible fluid is the exclusion process. Let $p(\cdot)$ be a finite range symmetric transition probability function on $\Z$. The simple symmetric exclusion process (SSEP) on $\mathbb{Z}$ associated to $p(\cdot)$ is a Feller process with configuration space $\Omega = \{0,1\}^{\mathbb{Z}}$ whose evolution can be described in the following way: initially each site of $\mathbb{Z}$ can be occupied or not by a particle, each particle of the system waits, independently of any other particle, an exponencial time of parameter one and at that time it choses another site according to $p(\cdot)$ and jumps to the chosen site if it is unoccupied. It is well-known that the SSEP has a hydrodynamic behavior under diffusive space-time scaling whose hydrodynamical equation is Laplace's equation, see \cite{kipnislandim}. 

\begin{figure}[htpb]
\label{fig:EP}
\includegraphics[viewport= 50 250 1200 550, width=100mm,height=30mm, clip]{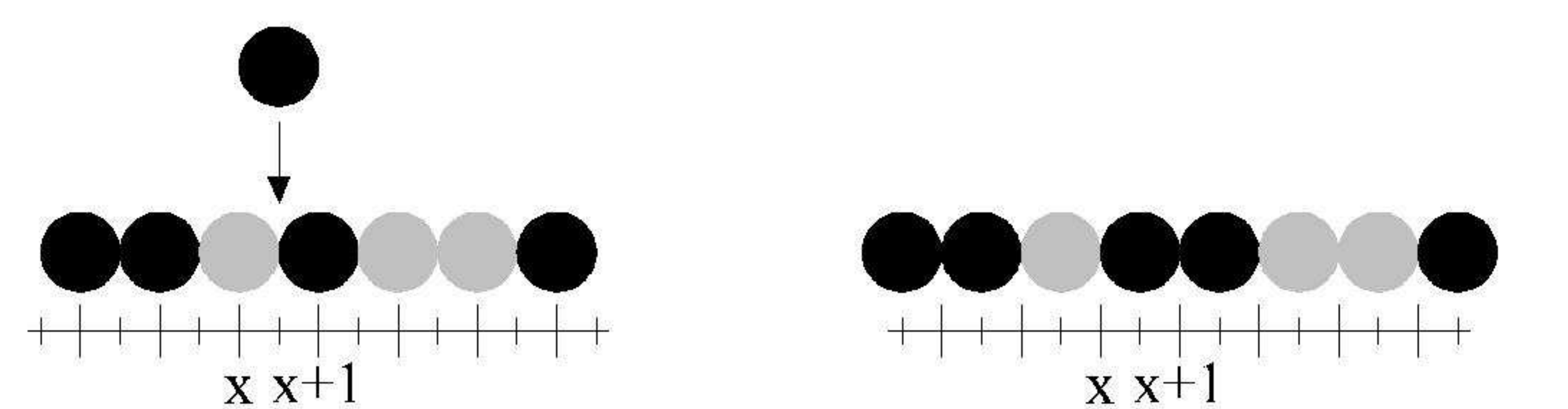}
\caption{A configuration for the exclusion process berore and after one extra black particle has entered}
\end{figure}

Consider the following alternative description of the exclusion process: each site is occupied by two types of particles, black and gray particles as in figure \ref{fig:EP}, each black and gray pair of particles, independently of any other such pair, exchange positions after an exponencial time of parameter $p(x-y)$ where $x$ is the position of the black particle and $y$ is the position of the gray one. In this way, black particles are associated to ones and gray particles to zeros for a configuration in $\Omega$. 

We describe now a pertubation of the exclusion process where new black particles are allowed to enter the system. The idea is to push a black particle between two nearest neighbor particles triggering an spread of particles centered in the position of the new particle, see figure \ref{fig:EP}. But after particles are translated, they are located in $\Z + \frac{1}{2}$ instead of $\Z$. To deal with this behavior, we enlarge the configuration space to $\tilde{\Omega} = \N \times (\Gamma_1 \cup \Gamma_2)$, where 
$$
\Gamma_1 = \Big\{ \eta \in \{0,1\}^{\Z \cup (\Z + \frac{1}{2})} : \eta(x) =0 \ \textrm{forall} \ x \in \Z + \frac{1}{2} \Big\}
$$
and
$$
\Gamma_2 = \Big\{ \eta \in \{0,1\}^{\Z \cup (\Z + \frac{1}{2})} : \eta(x) =0 \ \textrm{forall} \ x \in \Z \Big\}
$$
for $\Z +\frac{1}{2} = \{x+\frac{1}{2} : x \in \Z \}$. The system alternates between a SSEP on $\Z$ and a SSEP on $\Z + \frac{1}{2}$, where the interchanges are a result of a superposed dynamics describing the spread mechanism that opens space for one new particle. This dynamics is described as follows: 

\begin{enumerate}
\item[(i)] Let $h:[0,T] \times \mathbb{R} \ra \mathbb{R}_+$ be a smooth bounded function such that there exist constants $C>0$ and $\beta > 0$ such that
\begin{equation}
\label{eq:exp}
h(t,u) \le C e^{- \beta |u|} ,
\end{equation}
for every $u\in \mathbb{R}$ and $t \ge 0$. In particular,
\begin{equation}
\label{eq:h}
\sup_{t \in [0,T]} C(t) < \oo \, ,
\end{equation}
where
$$
C(t) := \int_{-\oo}^{+\oo} h(t,u) du \, ,
$$
which is a bounded Lebesgue measurable function.
\item[(ii)] Fix $N>0$ as the scaling parameter and let $s$ denote a fixed macroscopic time. If the state of the system at macroscopic time $s$ is $(n,\eta) \in \tilde{\Omega}$ with $\eta \in \Gamma_1$ (resp. $\eta \in \Gamma_2$), then we have that independently at each site $x \in \Z$ (resp. $x \in \Z+\frac{1}{2}$) and at rate $h(s,x/N)$ the system changes to $(n+1,\eta^\prime)$ with $\eta^\prime \in \Gamma_2$ (resp. $\Gamma_1$).
\item[(iii)] The configuration $\eta^\prime$, mentioned in (ii), is obtained from $\eta$ in the following way: all the system strictly at the right of site $x$ ($x$ included) is translated by $1/2$ units to the right; all the system on its left hand side is translated by $1/2$ units to the left; we set $\eta^\prime (x) = 0$ and $\eta^\prime (x-\frac{1}{2}) = 1$. 
\end{enumerate}

We are only going to consider initial cofigurations on $\tilde{\Omega}$ of the kind $(1,\eta)$ for $\eta\in \Gamma_1$, and therefore the system is always in $\Gamma_1$ for $n$ odd and in $\Gamma_2$ for $n$ even. Thus as mentioned before during each interval of time where $n$ is odd (resp. even) the system behaves as a SSEP on $\Z$ (resp. $\Z + \frac{1}{2}$). From condition (\ref{eq:h}), we only have a finite increase of mass in finite time which is required for the system to be well defined. The system we obtain will be called the Exclusion Process with centered spreading mechanism (EPCS). 

\smallskip
In order to prove the hydrodynamic behavior of the EPCS we perform two transformations on the configuration space obtaining a system we can deal with more easily. The first transformation is obtained by reversing the roles between zeros and ones. The second transformation makes the process into a SSEP on $\mathbb{Z}$ with a superposed dynamics. To describe this dynamics set $N$ as the scaling parameter and $s$ as a fixed macroscopic time. Let $b:[0,T] \times \mathbb{R} \ra \mathbb{R}_+$ be related to $h(\cdot,\cdot)$ by the equality
\begin{equation}
\label{eq:bh}
b(t,u) = h \Big( t , u - \int_0^t \frac{C(s)}{2} ds \Big) .
\end{equation} 
Then for a site $x\in \mathbb{Z}$ at rate $b(s,x/N)$, independently of any other site, the particles at the right hand side of $x$ are translated one unit to the right. 

The aim of the first transformation is to obtain a system for which we do not need to deal with incoming particles. The second transformation is to simplify the configuration space obtaining a gradient modified exclusion process. The system we get will be called the exclusion process with right sided spreading mechanism (EPRS).

The intuitive reasoning behind (\ref{eq:bh}) is that the EPRS is obtained from the EPCS by translating the system one half units to the right each time mass increases by one unit and thus $b$ at time $t$ should be the translation of $h$ in the space variable backwards by a half of the total mass that has entered the system before time $t$. Futhermore 
\begin{equation}
\label{eq:b}
\int_{-\oo}^{+\oo} b(t,u) du = C(t) \, , 
\end{equation}
and from (\ref{eq:h}) and the formal description of the EPRS through its generator we can show that the process is well-defined, see Section \ref{sec:serp}.

\smallskip
The EPRS is a weak pertubation of the SSEP by shift operators and to prove its hydrodynamic behavior we have to deal with the particularities of this pertubation. The proof presented here is based on the entropy method of Guo, Papanicolaou and Varadhan \cite{guopapanicolaouvaradhan} and an adaptation by Kipnis, Olla and Varadhan \cite{kipnisollavaradhan} to weak pertubations of the SSEP, but in \cite{kipnisollavaradhan} the pertubation of the Kawazaki dynamics is by Glauber instead of the shift operators for the EPRS. Hence this paper adds to the efforts to extend the results on hydrodynamics to new classes of non-conservative systems. Concerning the proof of the hydrodynamic behavior of non-conservative interacting particle systems, we have that similar arguments have also been employed in \cite{valle}. Other papers as \cite{landimollavolchan} and \cite{landimvalle} also deal with adaptations of the entropy method and in \cite{mourragui} the relative entropy method is applied. 

The EPCS and its hydrodynamic behavior are described in Section \ref{sec:hrp} and the EPRS and its hydrodynamic behavior in Section \ref{sec:serp}. The last sections, Sections \ref{sec:proofRP} and \ref{sec:proof}, are devoted to the proofs of the hydrodynamic limits for the EPCS and the EPRS respectively.

\medskip
%%%%%%%%%%%%%%%%%%%%%%%%%%%%%%%%%%%%%%%%%%%%%%%%%%%%%%%%%%%%%%%%%%%%%%%%%%%%
\section{Hydrodynamics for the EPCS}
%%%%%%%%%%%%%%%%%%%%%%%%%%%%%%%%%%%%%%%%%%%%%%%%%%%%%%%%%%%%%%%%%%%%%%%%%%%%
\label{sec:hrp}

The EPCS is formally described through its generator which is an operator on the space of local real functions with domain $\tilde{\Omega}$. Denote it by
$$
\mathcal{G}^N_s = N^2 \tilde{L} + \tilde{L}_r^{N,s} \, .
$$
The operator $\tilde{L}$ is related to the SSEP dynamics associated to a finite range symmetric transition probability $p(\cdot)$ and is equal to  
$$
\textbf{1}_{\eta \in \Gamma_1} \sum_{x \in \mathbb{Z}} \sum_{z \in \mathbb{Z}} p(z) \{ L_{x,x+z} + L_{x+z,x} \} + \textbf{1}_{\eta \in \Gamma_2} \sum_{x \in \mathbb{Z}+\frac{1}{2}} \sum_{z \in \mathbb{Z}} p(z) \{ L_{x,x+z} + L_{x+z,x} \}
$$
where, for every local function $F:\tilde{\Omega} \ra \mathbb{R}$ we have 
$$
L_{x,y} F(n,\eta) = \eta(x) [1-\eta(y)] \ [F(n,\eta^{x,y}) - F(n,\eta)],
$$
for every $(n,\eta) \in \tilde{\Omega}$ and $x,y \ge 1$, with $\eta^{x,y}$ being obtained from configuration $\eta$ with spins at $x$ and $y$ interchanged:
$$
\eta^{x,y}(z)= \left\{ \begin{array}{l}
\eta(y) \, , \ \textrm{if } z=x \\
\eta(x) \, , \ \textrm{if } z=y \\
\eta(z) \, , \ {\textrm{otherwise}} \, . 
\end{array} \right. 
$$
The operator $\tilde{L}_r^{N,s}$ is related to the spread of the mass around the position where a raindrop falls: 
$$
\tilde{L}_r^{N,s} = \textbf{1}_{\eta \in \Gamma_1} \sum_{x \in \mathbb{Z}} \tilde{L}_{r,x}^{N,s} +  \textbf{1}_{\eta \in \Gamma_2} \sum_{x \in \mathbb{Z}+\frac{1}{2}} \tilde{L}_{r,x}^{N,s}  
$$
where, for every local function $F:\tilde{\Omega} \ra \mathbb{R}$
$$
\tilde{L}_{r,x}^{N,s} F (n,\eta) = h(s,x/N) [F(n+1,\tilde{\tau}_x \eta)-F(n,\eta)] ,
$$
where, for every $y \in \Z \cup (\Z+\frac{1}{2})$,
$$
(\tilde{\tau}_x \eta)(y) = \left\{ \begin{array}{cl}
\eta(y-\frac{1}{2}) & \!\!\!\!\! , \ \textrm{if } y>x \\
\eta(y+\frac{1}{2}) & \!\!\!\!\! , \ \textrm{if } y<x-\frac{1}{2} \\
1 & \!\!\!\!\! , \ \textrm{if } y=x-\frac{1}{2} \\
0 & \!\!\!\!\! , \ \textrm{if } y=x \, .
\end{array} \right. 
$$

Note that condition (\ref{eq:b}) implies, as an application of the Borel-Cantelli Lemma that in finite time we only have a finite increase of mass. Thus, from the results presented in \cite{liggett} there exists a non-homogeneous Feller process associated to the generator $\mathcal{G}^N_s$ which we are going to denote by $(n^N_s,\eta^N_s)_{s\ge 0}$.

\medskip
%%%%%%%%%%%%%%%%%%%%%%%%%%%%%%%%%%%%%%%%%%%%%%%%%%%%%%%%%%%%%%%%%%%%%%%%%%%%
\subsection{Weak solutions of the EPCS hydrodynamic equation}
%%%%%%%%%%%%%%%%%%%%%%%%%%%%%%%%%%%%%%%%%%%%%%%%%%%%%%%%%%%%%%%%%%%%%%%%%%%%
\label{sec:weaksolutionsRP}

The Hydrodynamic equation for the EPCS is the convective diffusion equation
\begin{equation}
\label{eq:eqhydRP}
\left\{ 
\begin{array}{l}
\partial_t \rho(t,u) = \sigma^2 \Delta \rho(t,u) - \frac{1}{2} \partial_u  (\gamma(t,u) \rho(t,u)) + h(t,u), \ (t,u)\in (0,T) \times \mathbb{R}, \\ 
\rho(0,u) = \rho_0(u), \ u \in \mathbb{R} , 
\end{array}
\right.
\end{equation}
where
$$
\gamma(t,u) = \int_{-\oo}^u h(s,v) dv - \int_{u}^{+\oo} h(s,v) dv  \, , \quad \sigma^2 = \frac{1}{2} \sum_{z\in \mathbb{Z}} z^2 p(z) \, .
$$

In order to remove the external field $h(t,u)$ from equation (\ref{eq:eqhydRP}) we apply the transformation $\tilde{\rho} = 1 - \rho$. This is also required to obtain the hydrodynamical behavior of the EPCS from that of the EPRS, since from the microscopic point of view it is equivalent to reverse the roles between zeros and ones. Thus, we say that a bounded measurable function $\rho:[0,T) \times \mathbb{R} \ra \mathbb{R}$ is said to be a weak solution of equation (\ref{eq:eqhydRP}) if $\tilde{\rho}:[0,T) \times \mathbb{R} \ra \mathbb{R}$ is a weak solution of 
\begin{equation}
\label{eq:eqhydRP2}
\left\{ 
\begin{array}{l}
\partial_t \tilde{\rho}(t,u) = \sigma^2 \Delta \tilde{\rho}(t,u) - \frac{1}{2} \partial_u  (\gamma(t,u) \tilde{\rho}(t,u)), \ (t,u)\in (0,t) \times \mathbb{R}, \\ 
\tilde{\rho}(0,u) = \tilde{\rho}_0(u), \ u \in \mathbb{R} . 
\end{array}
\right.
\end{equation}

Denote by $C^{0,1}_0 ([0,T] \times \mathbb{R})$ the space of real continuous functions with compact support on $[0,T] \times \mathbb{R}$ which are continuously differentiable in the second variable and by $C^{1,2}_0 ([0,T] \times \mathbb{R})$ the space of real functions with compact support on $[0,T] \times \mathbb{R}$ which are continuously differentiable in the first variable and twice continuously differentiable in the second variable. Let $\tilde{\rho}_0:\mathbb{R} \ra \mathbb{R}$ be a fixed bounded measurable function. We say that a bounded measurable function $\tilde{\rho}:[0,T) \times \mathbb{R} \ra \mathbb{R}$ is a weak solution of equation (\ref{eq:eqhydRP2}) if
\begin{itemize}
\item[\textbf{(a)}] $\tilde{\rho}(t,u)$ is absolutely continuous in the space variable and $\partial_u \tilde{\rho}(t,u)$ is a locally square integrable function on $(0,T) \times \mathbb{R}$ such that for all $0\le t \le T$ and for every function $G \in C^{0,1}_0 ([0,T] \times \mathbb{R})$ we have that
$$
\int_0^T ds \int_{\mathbb{R}} du \, G(s,u) \partial_u \tilde{\rho}(s,u) = 
- \int_0^T ds \int_{\mathbb{R}} du \, \partial_uG(s,u) \tilde{\rho}(s,u) \, .
$$
\item[\textbf{(b)}] For every function $G \in C^{1,2}_0 ([0,T] \times \mathbb{R})$ and every $t \in [0,T]$ we have that
\begin{eqnarray}
\lefteqn{\qquad \int_{\mathbb{R}} du \, G(t,u) \tilde{\rho}(t,u) - \int_{\mathbb{R}} du \, G(0,u) \tilde{\rho}_0(u) =
\int_0^t ds \int_{\mathbb{R}} du \, \partial_s G(s,u) \tilde{\rho}(s,u) } \nn \\
& & \qquad + \ \int_0^t ds \left\{  -\sigma^2 \int_{\mathbb{R}} du \, \partial_u G(s,u) \partial_u \tilde{\rho}(s,u) + \frac{1}{2} \int_{\mathbb{R}} du \, \partial_u G(s,u) \gamma(s,u) \tilde{\rho}(s,u) \right\} . \nn
\end{eqnarray}  
\end{itemize}

Existence and regularity for equation (\ref{eq:eqhydRP2}), and thus for equation (\ref{eq:eqhydRP}), follows from existence and regularity for equation (\ref{eq:eqhyd1}), this is the contend of Lemma \ref{lemma:RPtoMRP}. The required uniqueness result is given by theorem 5.2 in \cite{lsu}.

\medskip
%%%%%%%%%%%%%%%%%%%%%%%%%%%%%%%%%%%%%%%%%%%%%%%%%%%%%%%%%%%%%%%%%%%%%%%%%%%%%
\subsection{The EPCS hydrodynamic behavior}
%%%%%%%%%%%%%%%%%%%%%%%%%%%%%%%%%%%%%%%%%%%%%%%%%%%%%%%%%%%%%%%%%%%%%%%%%%%%%
\label{sec:hbRP}
We denote by $\mathcal{P}(\Omega)$ the set of probability measures on $\Omega$. Let $\nu_\alpha \in \mathcal{P}(\Omega)$ be the Bernoulli product measure of parameter $\alpha \in [0,1]$ on $\Omega$. Given any two probabilities $\mu$, $\nu$ in $\mathcal{P}(\Omega)$, we denote by $H(\mu|\nu)$ the relative entropy of $\mu$ with respect to $\nu$:
$$
H(\mu|\nu)=\sup_f \left\{ \int fd\mu - \log \int e^f d\nu \right\} ,
$$
where the supremum is carried over all bounded continuous real functions. 

A sequence of measures $(\mu_N)_{N \ge 1}$ in $\mathcal{P}(\Omega)$ is associated to a initial profile $\rho_0 : \mathbb{R} \ra \mathbb{R}$ if for every $\delta > 0$ and every continuous function $G:\mathbb{R} \ra \mathbb{R}$ with compact support
$$
\lim_{N\ra +\oo}
\mu^N \left( \left| \frac{1}{N} \sum_{z \in \mathbb{Z}} G(z/N) \eta(z) - \int du \, G(u) \rho_0 (u) \right| > \delta \right) = 0 \, .
$$

\medskip
Let $D(\mathbb{R}_+,\tilde{\Omega})$ denote the space of right continuous functions with left limits on $\tilde{\Omega}$ endowed with the Skorohod topology. For each probability measure $m$ on $\tilde{\Omega}$ and $N\ge 1$, denote by $\tilde{\mathbb{P}}_{m}^N$ the probability measure on $D(\mathbb{R}_+,\tilde{\Omega})$ induced by the Markov process $(n^N_s,\eta^N_s)_{s\ge 0}$ with generator $\mathcal{G}^N_s$ and with initial measure $m$.

For a measure $\mu$ in $\mathcal{P}(\Omega)$, let $\mathcal{T}(\mu)$ be the measure on $\tilde{\Omega}$ induced by $\mu$ through the transformation $\xi \in \Omega \ra (1,\eta) \in \tilde{\Omega}$ where $\eta$ is the configuration in $\Gamma_1$ which is equal to $\xi$ on $\mathbb{Z}$.

\smallskip
The hydrodynamical behavior of the EPCS is given by the following result:

\begin{theorem} %%%%%%%%%%%%%%%%%%%%%%%%%%%%%%%%%%%%%%%%%%%%%%%%%%%%%%%%%
\label{theorem:hbRP}
Fix a sequence of initial probability measures $(\mu_N)_{N \ge 1}$ in $\mathcal{P}(\Omega)$ associated to an initial profile $\rho_0: \mathbb{R} \ra \mathbb{R}_+$, bounded above by $1$, such that $H(\mu^N|\nu_\alpha) \le CN$ for some $\alpha \in (0,1)$. Then, for any continuous function $G: \mathbb{R} \ra \mathbb{R}$ with compact support, any $\delta>0$ and $0<t<T$ 
$$
\lim_{N\ra {\oo}} \tilde{\mathbb{P}}_{\mathcal{T}(\mu^{\scriptscriptstyle{N}})}^N \left[ \left| \frac{1}{N} \sum_{z \in \Z \cup (\Z + \frac{1}{2})} G(z/N) \eta_t(z) - \int duG(u)\rho(t,u) \right| \ge \delta \right] = 0
$$
where $\rho$ is the unique weak solution of the convective diffusion equation (\ref{eq:eqhydRP}).
\end{theorem} %%%%%%%%%%%%%%%%%%%%%%%%%%%%%%%%%%%%%%%%%%%%%%%%%%%%%%%%%%%% 

\medskip
%%%%%%%%%%%%%%%%%%%%%%%%%%%%%%%%%%%%%%%%%%%%%%%%%%%%%%%%%%%%%%%%%%%%%%%%%%%%
\section{Hydrodynamics for the EPRS}
%%%%%%%%%%%%%%%%%%%%%%%%%%%%%%%%%%%%%%%%%%%%%%%%%%%%%%%%%%%%%%%%%%%%%%%%%%%%
\label{sec:serp}

The EPRS is formally described through its generator which is an operator on the space of local real functions with domain $\Omega$. Denote it by 
$$
\mathcal{L}^N_s = N^2 L + L_r^{N,s}.
$$ 
Here $L$ is related to the motion of particles in a SSEP associated to a finite range symmetric transition probability $p(\cdot)$: 
$$
L = \sum_{x \in \mathbb{Z}} \sum_{z \in \mathbb{Z}} p(z) \{ L_{x,x+z} + L_{x+z,x} \}
$$
where, for every local function $F:\Omega \ra \mathbb{R}$ and every $x,y \ge 1$,
$$
L_{x,y} F(\xi) = \xi(x) [1-\xi(y)] \ [F(\xi^{x,y}) - F(\xi)],
$$
and $\xi^{x,y}$ is the configuration with spins at $x$ and $y$ interchanged:
$$
\xi^{x,y}(z)= \left\{ \begin{array}{l}
\xi(y) \, , \ \textrm{if } z=x \\
\xi(x) \, , \ \textrm{if } z=y \\
\xi(z) \, , \ {\textrm{otherwise}} \, . 
\end{array} \right. 
$$
The operator $L_r^{N,s}$ is related to the translations of the system at the right of a given site:
$$
L_r^{N,s} = \sum_{x \in \mathbb{Z}} L_{r,x}^{N,s} 
$$
where, for every local function $F:\Omega \ra \mathbb{R}$
$$
L_{r,x}^{N,s} F (\xi) = b(s,x/N) [F(\tau_x \xi)-F(\xi)] ,
$$
with 
$$
(\tau_z \xi)(x) = \left\{ \begin{array}{cl}
\xi(x-1) & \!\!\!\!\! , \ \textrm{if } x>z \, , \\
0 & \!\!\!\!\! , \ \textrm{if } x=z \, , \\
\xi(x) &  \!\!\!\!\! , \ \textrm{if } x<z \, .
\end{array} \right. 
$$
We refer to Liggett's book \cite{liggett}, to a proof that $\mathcal{L}^N_s$ is the generator of a Feller semigroup. Therefore, there exists a non-homogeneous Feller Process associated to the generator $\mathcal{L}^N_s$ which we are going to denote by $(\xi^N_s)_{s \ge 0}$.

\medskip
%%%%%%%%%%%%%%%%%%%%%%%%%%%%%%%%%%%%%%%%%%%%%%%%%%%%%%%%%%%%%%%%%%%%%%%%%%%%
\subsection{Weak solutions of the EPRS hydrodynamic equation}
%%%%%%%%%%%%%%%%%%%%%%%%%%%%%%%%%%%%%%%%%%%%%%%%%%%%%%%%%%%%%%%%%%%%%%%%%%%%
\label{sec:weaksolutions}

The hydrodynamic equation of the EPRS is also a convective diffusion equation:
\begin{equation}
\label{eq:eqhyd1}
\left\{ 
\begin{array}{l}
\partial_t \zeta(t,u) = \sigma^2 \Delta \zeta(t,u) - \partial_u  (a(t,u) \zeta(t,u)), \ (t,u)\in \mathbb{R}_+ \times \mathbb{R}, \\ 
\zeta(0,u) = \zeta_0(u), \ u \in \mathbb{R} , 
\end{array}
\right.
\end{equation}
where
$$
a(t,x) = \int_{-\oo}^x b(s,u) du\, , \quad \sigma^2 = \frac{1}{2} \sum_{z\in \mathbb{Z}} z^2 p(z) \, .
$$

For a fixed bounded mesurable function $\zeta_0:\mathbb{R} \ra \mathbb{R}$, a bounded measurable function $\zeta:[0,T) \times \mathbb{R} \ra \mathbb{R}$ is said to be a weak solution of equation (\ref{eq:eqhyd1}) if
\begin{itemize}
\item[\textbf{(a)}] $\zeta(t,u)$ is absolutely continuous in the space variable and $\partial_u\zeta(t,u)$ is a locally square integrable function on $(0,T) \times \mathbb{R}$ such that for all $0\le t \le T$ and for every function $G \in C^{0,1}_0 ([0,T] \times \mathbb{R})$ we have that
$$
\int_0^T ds \int_{\mathbb{R}} du \, G(s,u) \partial_u\zeta(s,u) = 
- \int_0^T ds \int_{\mathbb{R}} du \, \partial_uG(s,u) \zeta(s,u) \, .
$$
\item[\textbf{(b)}] For every function $G \in C^{1,2}_0 ([0,T] \times \mathbb{R})$ and every $t \in [0,T]$,
\begin{eqnarray}
\lefteqn{\qquad \int_{\mathbb{R}} du \, G(t,u)\zeta(t,u) - \int_{\mathbb{R}} du \, G(0,u)\zeta_0(u) =
\int_0^t ds \int_{\mathbb{R}} du \, \partial_s G(s,u) \zeta(s,u) } \nn \\
& & \qquad + \ \int_0^t ds \left\{  -\sigma^2 \int_{\mathbb{R}} du \, \partial_u G(s,u) \partial_u \zeta(s,u) + \int_{\mathbb{R}} du \, \partial_u G(s,u) a(s,u) \zeta(s,u) \right\} . \nn
\end{eqnarray}  
\end{itemize}

Existence and regularity for equation (\ref{eq:eqhyd1}) follows from the proof of the hydrodynamical behavior of the system described in Theorem \ref{theorem:hbep} on Section \ref{sec:hbep}. The required uniqueness result is given by theorem 5.2 in \cite{lsu}.

\medskip
%%%%%%%%%%%%%%%%%%%%%%%%%%%%%%%%%%%%%%%%%%%%%%%%%%%%%%%%%%%%%%%%%%%%%%%%%%%%%
\subsection{The EPRS hydrodynamic behavior}
%%%%%%%%%%%%%%%%%%%%%%%%%%%%%%%%%%%%%%%%%%%%%%%%%%%%%%%%%%%%%%%%%%%%%%%%%%%%%
\label{sec:hbep}
Let $D(\mathbb{R}_+,\Omega)$ denote the space of right continuous functions with left limits on $\Omega$ endowed with the Skorohod topology. For each probability measure $\mu$ on $\Omega$ and $N\ge 1$, denote by $\mathbb{P}_{\mu}^N$ the probability measure on $D(\mathbb{R}_+,\Omega)$ induced by the Markov process $\xi_t^N$ with generator $\mathcal{L}^N_s$ and with initial measure $\mu$. Recall also the definitions of $\nu_\alpha$, the relative entropy from Section \ref{sec:hbRP} and of a family of probability measures on $\Omega$ associated to a density profile.

The hydrodynamical behavior of the EPRS is given by the following result:

\begin{theorem} %%%%%%%%%%%%%%%%%%%%%%%%%%%%%%%%%%%%%%%%%%%%%%%%%%%%%%%%%
\label{theorem:hbep}
Fix a sequence of initial probability measures $(\mu^N)_{N\ge 1}$ on $\Omega$ associated to a initial profile $\zeta_0: \mathbb{R} \ra \mathbb{R}$, bounded above by $1$, such that $H(\mu^N|\nu_\alpha) \le CN$ for some $\alpha \in (0,1)$. Then, for any continuous function $G: \mathbb{R} \ra \mathbb{R}$ with compact support, any $\delta>0$ and $0<t<T$ 
$$
\lim_{N\ra {\oo}} \mathbb{P}_{\mu^{\scriptscriptstyle{N}}}^N \left[ \left| \frac{1}{N} \sum_{z \in \mathbb{Z}} G(z/N) \xi_t(z) - \int duG(u)\zeta(t,u) \right| \ge \delta \right] = 0
$$
where $\zeta$ is the unique weak solution of (\ref{eq:eqhyd1}).
\end{theorem} %%%%%%%%%%%%%%%%%%%%%%%%%%%%%%%%%%%%%%%%%%%%%%%%%%%%%%%%%%%% 

\medskip
%%%%%%%%%%%%%%%%%%%%%%%%%%%%%%%%%%%%%%%%%%%%%%%%%%%%%%%%%%%%%%%%%%%%%%%%%%%%%
\section{From the Hydrodynamic for the EPRS to the hydrodynamic for the EPCS}
%%%%%%%%%%%%%%%%%%%%%%%%%%%%%%%%%%%%%%%%%%%%%%%%%%%%%%%%%%%%%%%%%%%%%%%
\label{sec:proofRP}
\setcounter{equation}{0}

This Section is devoted to the proof of Theorem \ref{theorem:hbRP}. The proof is based on a coupling and thus comparative arguments between the EPCS and the EPRS, therefore it follows from Theorem \ref{theorem:hbep}, which is shown in Section \ref{sec:proof}. For the coupling we use an auxiliary process which is obtained by translating to the left the EPRS by a half times the number of translations of the system that have already ocurred and then by reversing the roles between zeros and ones. In this way, if we put in correspondence the number of particles entering the EPCS with the translations in the EPRS, we have that the auxiliary process is basicaly the EPCS except for corrections in the positions of particles that are negligigle in probability as $N\ra +\oo$, where $N$ is the scaling parameter.

Recall from the statement of Theorem \ref{theorem:hbRP} that we will be considering a family of probability measures $(\mu^N)_{N\ge 1}$ on $\Omega$ satisfying some conditions imposed there. Also from Section \ref{sec:hbRP}, this family can be regarded as a family of probability measures on $\tilde{\Omega}$ through the correspondence between $\Gamma_1$ and $\Omega$. Thus it represents the initial distribuition for both the EPRS and the EPCS in this section. 

Aiming at the strategy described above, let us first make some considerations about the number of translations of the EPRS. At scaling $N$ and for each $x \in \Z$, we denote by $(W^{x,N}_t : 0 \le t \le T)$ the independent non-homogeneous Poisson processes with rate $b(s,x/N)$, $0\le s \le T$, which represents the number of translations to the right of the EPRS ocurring from site $x$ before time $t$.
Next lemma relates in a weak sense the total number of translations of the EPRS with the macroscopic rate at which the system is translated.

\begin{lemma}
\label{lemma:wlln}
For every $\delta > 0$ and $0 \le t \le T$ we have that
$$
\lim_{N\ra {\oo}} \mathbb{P}_{\mu^{\scriptscriptstyle{N}}}^N \left[ \left| \frac{1}{N} \sum_{x \in \Z} W^{x,N}_t - \int_0^t C(s) ds \right| \ge \delta \right] = 0\, .
$$
\end{lemma}

\no \textbf{Proof:} We show that for every $u \in \mathbb{R}$
$$
\lim_{N\ra {\oo}} \mathbb{P}_{\mu^{\scriptscriptstyle{N}}}^N \left[ \left| \frac{1}{N} \sum_{x \le [uN]} W^{x,N}_t - \int_0^t \int_{-\oo}^u h(s,v) \, dvds \right| \ge \delta \right] = 0 \, .
$$
It is clear in the proof that we may replace $u$ by $+\oo$ in the previous expression to get the limit in the statement. Since convergence in probability to a constant follows from convergence in distribution, we only have to show that $\frac{1}{N} \sum_{x \le [uN]} W^{x,N}_t$ converges in distribution to $\int_0^t \int_{-\oo}^u h(s,v) \, dvds$. We are going to estalish the convergence of the appropriate characteristic functions. Since, for a Poisson process $W$ with time varying rate $\lambda(s)$ we have that
$$
E[\exp\{i \gamma W_t \}] = \exp \Big\{ (e^{i\gamma} - 1) \int_0^t \lambda(s) ds \Big\},
$$
we get from independence that 
\begin{eqnarray*}
E\Big[ \exp \Big\{ \frac{i \gamma}{N} \sum_{x \le [uN]} W_t^{x,N} \Big\} \Big] & = & \Pi_{x \le [uN]} E \Big[ \exp \Big\{ \frac{i \gamma \, W^{x,N}_t}{N} \Big\} \Big] \\
& = & \Pi_{x \le [uN]} \exp \Big\{ (e^{\frac{i \gamma}{N}} - 1) \int_0^t h(s,x/N) ds \Big\} \\
& = & \exp \Big\{ \sum_{x \le [uN]} (e^{\frac{i \gamma}{N}} - 1) \int_0^t h(s,x/N) ds \Big\},
\end{eqnarray*}
for every $\gamma \in \mathbb{R}$. Now
$$
\sum_{x \le [uN]} (e^{\frac{i \gamma}{N}} - 1) \int_0^t h(s,x/N) ds \, \ra \, i \gamma \int_0^t \int_{-\oo}^u h(s,v) dvds \, , \textrm{ as } N \ra \oo \, , 
$$
and then
$$
E\Big[ \exp \Big\{ \frac{i \gamma}{N} \sum_{x \le [uN]} W_s^{x,N} \Big\} \Big] \ra \exp \Big\{ i \gamma \int_0^t \int_{-\oo}^u h(s,u) duds \Big\} \, , \textrm{ as } N \ra \oo \, ,
$$
for every $\gamma > 0$. Therefore we have proved the required convergence.
$\square$

\bigskip

Let $(\xi^N_s)_{s\ge 0}$ denote the EPRS. To obtain the auxiliary system we consider the transformation from $D(\mathbb{R}_+,\Omega)$ to $D(\mathbb{R}_+,\tilde{\Omega})$ which associates to the configuration $\xi^N_s$ the configuration $(\hat{n}^N_s,\hat{\eta}^N_s)$ given by
$$
\hat{n}^N_s = \sum_{z \in \Z} W^{z,N}_s
$$
and
$$
\hat{\eta}^N_s(x) = 
\left\{ \begin{array}{ll}
1 - \xi^N_s \Big( x + \frac{\hat{n}^N_s}{2} \Big)  \!\!\!\!\!\! &, \textrm{ if } \big( x + \frac{\hat{n}^N_s}{2} \big) \in \Z \\
0 \!\!\!\!\!\! &,  \textrm{ otherwise.}
\end{array} \right.
$$
Thus
$$
\frac{1}{N} \!\! \sum_{x \in \Z \cup (\Z + \frac{1}{2})} \!\!\!\! G(x/N) \hat{\eta}^N_s(x) 
$$
is equal to
\begin{eqnarray}
\label{eq:etaprime}
\lefteqn{\frac{1}{N} \sum_{x \in \Z } G\left(\frac{1}{N} \Big( x + \frac{\hat{n}^N_s}{2} - \left\lfloor \frac{\hat{n}^N_s}{2} \right\rfloor \Big) \right) \Big[ 1 - \xi^N_s \Big( x + \left\lfloor \frac{\hat{n}^N_s}{2} \right\rfloor  \Big) \Big] } \nn \\
\!\! & = & \!\! \frac{1}{N} \sum_{x \in \Z} G\Big( \frac{1}{N} \Big( x + \frac{\hat{n}^N_s}{2} - 2 \left\lfloor \frac{\hat{n}^N_s}{2} \right\rfloor \Big) \Big) [1-\xi^N_s(x)] \, .
\end{eqnarray}
By Theorem \ref{theorem:hbep} and Lemma \ref{lemma:wlln} the term in the rightmost side of the previous equation converges in probability to
\begin{equation}
\label{eq:etaprime1}
\int_{\mathbb{R}} G \Big( u - \int_0^s \frac{C(r)}{2} dr \Big) (1-\zeta(s,u)) du 
\end{equation}
for all $0\le s \le T$, where $\zeta:[0,T)\times \mathbb{R} \ra \mathbb{R}$ is the weak solution of (\ref{eq:eqhyd1}). We want to show that the same holds for the EPCS replacing the auxiliary process. To show this assertion we couple the EPCS with the auxiliary process in a proper way. Our coupling will allow first and second class particles on both systems. The difference between them is that first class particles have priority over second class particles this means that if a first class particle attempts to jump over a second class particle then they exchange positions. On the other hand, if a second class particle attempts to jump over a first class particle the jump is suppressed. The description of the coupling follows below:
\begin{enumerate}
\item[(i)] Both systems start with the same configuration.  For both systems, particles in the initial configuration are labeled in an ordered way from left to right. All these particles are first class, indeed only the first class particles are going to be labeled.
\item[(ii)] A jump attempt of a single first class particle occurs for the EPCS if and only if a jump attempt in the same direction occurs for the first class particle with the same label in the auxiliary process. Second class particles on both systems move independently of any other particle.
\item[(iii)] The appearance of a particle and the induced translation may occur in only one of the systems at a time or in both of them simultaneously. Recall that for the EPCS a particle appears at site $x-\frac{1}{2}$ at time $s$ with rate $h(s,x/N) = b(s,\frac{x}{N}+ \int_0^s \frac{C(r)}{2} dr)$ and for the auxiliary process with rate $b(s,\frac{x}{N}+ \frac{\hat{n}^N_s}{2N})$. In order to deal with the fact that a coordinate process might be evolving on $\Z$ while the other is evolving on $(\Z + \frac{1}{2})$, define $x^\prime = x$ if both coordinate processes have received the same amount of new particles and are evolving on the same lattice, and $x^\prime = x + \frac{1}{2}$ otherwise. We have one of three possibilities:
\begin{enumerate}
\item[(1)] If $h(s,x/N)>b(s,\frac{x^\prime}{N}+ \frac{\hat{n}^N_s}{2N})$, only the EPCS is translated around site $x$ at time $s$ at random rate $h(s,x/N) - b(s,\frac{x^\prime}{N}+ \frac{\hat{n}^N_s}{2N})$ with a particle appearing at site $x-\frac{1}{2}$. This particle becomes a second class particle for the EPCS.
\item[(2)] If $h(s,x/N)<b(s,\frac{x^\prime}{N}+ \frac{\hat{n}^N_s}{2N})$, only the auxiliary process is translated around site $x^\prime$ at time $s$ with a particle appearing at site $x^\prime-\frac{1}{2}$ at random rate $b(s,\frac{x^\prime}{N}+\frac{\hat{n}^N_s}{2N}) - h(s,x/N)$. This particle becomes a second class particle for the auxiliary process.
\item[(3)] Both systems are translated simultaneously at time $s$, the EPCS around site $x$ with a particle appearing at site $x-\frac{1}{2}$ and the EPRS around site $x^\prime$ with a particle appearing at site $x^\prime - \frac{1}{2}$. This happens at random rate $\min \{h(s,x/N),b(s,\frac{x^\prime}{N}+ \frac{\hat{n}^N_s}{2N})\}$. These particles are first class. Since we have new first class particles, we have to relabel them in a way to keep the particles ordered from left to right. Let us explain how we relabel these particles: before the transition, due to (1) and (2), we may have pairs of equally labeled particles that do not share the same positions. Let the positions of one of these pairs be denoted by $(z,w)$ with the first particle belonging to the EPCS and the other belonging to the auxiliary process. If $z<x$ and $w<x^\prime$ or $z\ge x$ and $w\ge x^\prime$ then the two particles remain to share the same label after transition. Now suppose that there exists a number $l$ (necessarily finite) of such pairs with positions $(z_j,w_j)$, $1\le j \le l$, such that $z_i<z_j<x$ and $x^\prime \le w_i < w_j$ or $x \le z_i < z_j$ and $w_i<w_j<x^\prime$, for every $i<j$. Now we consider new pairs $(z_1-\frac{1}{2},x^\prime-\frac{1}{2})$, $(z_2-\frac{1}{2},w_1+\frac{1}{2})$, ... , $(z_l-\frac{1}{2},w_{l-1}+\frac{1}{2})$, $(x-\frac{1}{2},w_l+\frac{1}{2})$ in the case $z_i<z_j<x$ and $x^\prime \le w_i < w_j$ or new pairs $(x-\frac{1}{2},w_1-\frac{1}{2})$, $(z_1+\frac{1}{2},w_2-\frac{1}{2})$, ... , $(z_{l-1}+\frac{1}{2},w_{l}-\frac{1}{2})$, $(z_l+\frac{1}{2},x^\prime-\frac{1}{2})$ in the case $x \le z_i < z_j$ and $w_i<w_j<x^\prime$. Finally give the same label to both particles in each of these pairs of positions. Note that in this transition the distance between two equally labeled particles can only decrease.
\end{enumerate}
\end{enumerate}
Following this we obtain a Markov process on $\tilde{\Omega} \times \tilde{\Omega}$ which will be represented by $((n^N_s,\eta^N_s),(\hat{n}_s^N,\hat{\eta}_s^N))_{s\ge 0}$ where $(n^N_s,\eta^N_s)$ is distribuited as the EPCS and $(\hat{n}^N_s,\hat{\eta}_s^N)$ as the auxiliary process, this will be called the coupled system. For a formal description, one can write explicitly the generator of the coupled system or use Harris graphical construction which are already standard in the literature, we refer here to Liggett's book \cite{liggett}.

Note that the maximum distance at time $s$ between two equally labeled particles in each one of the two coordinate processes of the coupled system is bounded by the total number of second class particles. Let this random number be denoted by $J^N_s$. We have that $(J^N_s)_{s\ge 1}$ is pure birth process with random rate
$$
\sum_{x\in \Z} \left| b\Big(s,\frac{x}{N}+ \int_0^s \frac{C(r)}{2} dr\Big) - b\Big(s,\frac{x^\prime}{N}+ \frac{\hat{n}^N_s}{2N}\Big) \right|.
$$ 
By Lemma \ref{lemma:wlln} and the exponencial bound on $b$ given in (\ref{eq:exp}), the integral in the time interval $[0,T]$ of the previous expression divided by $N$ goes to zero in probability as $N \ra \oo$. Therefore, for every $\epsilon>0$ and $0<\gamma<1$ we can fix $N_0=N_0(\epsilon)$ such that for $N>N_0$ the random variable $\frac{J^N_s}{N}$ is stochastically bounded on a set of probability greater than $\gamma$ by a Poisson distribuition with mean $\epsilon N$ and scaling parameter $N^{-1}$. This Poisson distribuition converges to zero in probability since it has characteristic function $\exp\{\epsilon N(e^{\frac{i}{N}}-1)\}$ which goes to $1$ as $N\ra \oo$ and $\epsilon \ra 0$. Hence $\frac{J^N_s}{N}$ also converges to zero in probability as $N \ra \oo$. Note that $|n^N_s - \hat{n}^N_s| \le J^N_s$, and therefore $\frac{|n^N_s - \hat{n}^N_s|}{N}$ also converges to zero in probability as $N \ra \oo$.

Now as in the description of the coupling we label the first class particles of the EPCS in a ordered way from left to right by denoting their positions at time $s$ by $X^N_k(s)$, for $k\in \Z$, and we do the same for the auxiliary process denoting the positions by $\hat{X}^N_k(s)$. We do the same with the second class particles denoting their positions respectively $Y^N_k(s)$, $1\le k \le m$, and $\hat{Y}^N_k(s)$, $1\le k \le l$, where $m = m(N,s)$ is the number of second class particles in the EPCS and $l = l(N,s)$ is the number of second class particles in the auxiliary process. In particular, $J^N_s = m+l$.

Finally we are in condition to compare integrals with respect to the empirical measures associated to $\eta$ and $\hat{\eta}$. We have that
\begin{equation}
\label{eq:competa}
\left| \frac{1}{N} \!\! \sum_{x \in \Z \cup (\Z + \frac{1}{2})} \!\!\!\! G(x/N) \eta^N_s(x) - \frac{1}{N} \!\! \sum_{x \in \Z \cup (\Z + \frac{1}{2})} \!\!\!\! G(x/N) \hat{\eta}^N_s(x) \right|
\end{equation}
is equal to
$$
\left| \frac{1}{N} \sum_{k \in \Z} \Big[ G\Big( \frac{X^N_k(s)}{N} \Big) -  G\Big( \frac{\hat{X}^N_k(s)}{N} \Big) \Big]+ \frac{1}{N} \sum_{k=1}^m G\Big( \frac{Y^N_k(s)}{N} \Big) -  \frac{1}{N} \sum_{k=1}^l G\Big( \frac{\hat{Y}^N_k(s)}{N} \Big) \right|.
$$
Since $G$ is a smooth function with compact support we have that this last expression is bounded above by  
$$
C(G) \, \Big( \sup_{k\in \Z} \frac{|X^N_k(s) - \hat{X}^N_k(s)|}{N} + \frac{J^N_s}{N} \Big) \le 2 \, C(G) \, \frac{J^N_s}{N} \, ,
$$
where $C(G)$ is some constant depending on $G$. Therefore (\ref{eq:competa}) converges to zero in probability as $N\ra +\oo$ and from (\ref{eq:etaprime}) and (\ref{eq:etaprime1})
any $\delta>0$ and $0<t<T$ 
$$
\lim_{N\ra {\oo}} \mathbb{P}_{\mathcal{T}(\mu^{\scriptscriptstyle{N}})}^N \left[ \left| \frac{1}{N} \sum_{z \in \Z \cup (\Z + \frac{1}{2})} G(z/N) \eta_t(z) - \int duG(u)\rho(t,u) \right| \ge \delta \right] = 0
$$
where $\rho$ is given by
\begin{equation} \label{eq:rhozeta}
\rho(t,u) = 1 - \zeta \Big( t , u + \int_0^t \frac{C(s)}{2} ds \Big), \textrm{ for every } (t,u) \in [0,T] \times \mathbb{R} \, .
\end{equation}
Now with the next lemma we end the proof:

\medskip
\begin{lemma}\label{lemma:RPtoMRP}Let $\zeta:[0,T) \times \mathbb{R} \ra \mathbb{R}$ be a weak solution of (\ref{eq:eqhyd1}). Then $\tilde{\rho}(\cdot,\cdot) = 1 - \rho(\cdot,\cdot)$, with $\rho(\cdot,\cdot)$ defined in (\ref{eq:rhozeta}), is a weak solution of (\ref{eq:eqhydRP2}).\end{lemma}

\no \textbf{Proof:} Note first that for a function $G \in C^{0,1}_0 ([0,T] \times \mathbb{R})$, the function $F$ defined by
\begin{equation} \label{eq:F}
F(t,u) = G \Big( t , u - \int_0^t \frac{C(s)}{2} ds \Big), \textrm{ for every } (t,u) \in [0,T] \times \mathbb{R} \, ,
\end{equation}
is also in $C^{0,1}_0 ([0,T] \times \mathbb{R})$. Thus we have that \textbf{(a)} in the definition of weak solutions of equation (\ref{eq:eqhydRP2}) follows for $\tilde{\rho}$ directly from (\ref{eq:rhozeta}) and \textbf{(a)} in definition of weak solutions of (\ref{eq:eqhyd1}) by a change of variables on both sides of the integral. 

\smallskip

Therefore it remains to verify that \textbf{(b)} in the definition of weak solutions of equation (\ref{eq:eqhydRP2}) follows from \textbf{(b)} in the definition of weak solutions of equation (\ref{eq:eqhyd1}). Fix a function $G \in C^{1,2}_0 ([0,T] \times \mathbb{R})$, we have that 
\begin{equation} \label{eq:G1}
\int_{\mathbb{R}} du \, G(t,u) \tilde{\rho}(t,u) - \int_{\mathbb{R}} du \, G(0,u) \tilde{\rho}_0(u) =
\int_{\mathbb{R}} du \, F (t,u) \zeta (t,u) - \int_{\mathbb{R}} du \, F (0,u) \zeta_0(u) \, ,
\end{equation}
where $F$ is as in (\ref{eq:F}). We can apply the equality \textbf{(b)} in the definition of weak solutions of (\ref{eq:eqhyd1}) if $F \in C^{1,2}_0 ([0,T] \times \mathbb{R})$, however we cannot guarantee that the first derivative in the first variable is differentiable since
$$
\partial_t F (t,u) = \partial_t G \Big( t , u - \int_0^t \frac{C(s)}{2} ds \Big) - \frac{C(s)}{2} \partial_u G\Big( t , u - \int_0^t \frac{C(s)}{2} ds \Big)
$$  
and $C$ is not necessarily continuous. To overcome this problem, we rely on some approximation results. Let $D:[0,T] \ra \mathbb{R}_+$ be the function $D(t) = \frac{1}{2} \int_0^t C(s) ds$. The function $D$ is a continuous increasing function with bounded, thus integrable, derivative. By Theorem 1 in Section 4.2 of \cite{evansgariepy}, there exists a family of smooth functions $D^\epsilon$, for $\epsilon>0$, obtained from $D$ by the convolution with standard mollifiers, such that $D^\epsilon$ converges uniformly to $D$  and $dD^\epsilon$ converges to $dD$ in $L^1([0,T])$. In accordance to (\ref{eq:F}), let $F^\epsilon$ be the smooth function $F^\epsilon (t,u) = G ( t , u - D^\epsilon ), \textrm{ for every } (t,u) \in [0,T] \times \mathbb{R}$. Since $F^\epsilon \in C^{1,2}_0 ([0,T] \times \mathbb{R})$ then  
\begin{eqnarray*}
\lefteqn{\qquad
\int_{\mathbb{R}} du \, F^\epsilon (t,u) \zeta (t,u) - \int_{\mathbb{R}} du \, F^\epsilon (0,u) \zeta_0(u) =
\int_0^t ds \int_{\mathbb{R}} du \, \partial_s F^\epsilon (s,u) \zeta(s,u) } \nn \\
& & \qquad + \ \int_0^t ds \left\{  -\sigma^2 \int_{\mathbb{R}} du \, \partial_u F^\epsilon (s,u) \partial_u \zeta(s,u) + \int_{\mathbb{R}} du \, \partial_u F^\epsilon (s,u) a(s,u) \zeta(s,u) \right\} . \nn
\end{eqnarray*}
Now, since $\zeta$ is bounded, $\partial_u \zeta$ is locally square integrable and from the definition of $F^\epsilon$, we can stablish by straightforward manipulations the convergence term to term as $\epsilon \ra 0$ in the right hand side of the previous equality. Therefore   
\begin{eqnarray*}
\lefteqn{\qquad
\int_{\mathbb{R}} du \, F (t,u) \zeta (t,u) - \int_{\mathbb{R}} du \, F (0,u) \zeta_0(u) =
\int_0^t ds \int_{\mathbb{R}} du \, \partial_s F (s,u) \zeta(s,u) } \nn \\
& & \qquad + \ \int_0^t ds \left\{  -\sigma^2 \int_{\mathbb{R}} du \, \partial_u F (s,u) \partial_u \zeta(s,u) + \int_{\mathbb{R}} du \, \partial_u F (s,u) a(s,u) \zeta(s,u) \right\} , \nn
\end{eqnarray*} 
where the right hand side can be rewritten as
\begin{eqnarray} \label{eq:F2}
\lefteqn{\!\!\!\!\!\!\!\!\!\!\!\!\!\!
\int_0^t ds \int_{\mathbb{R}} du \, \partial_s G(s,u) \tilde{\rho}(s,u) 
+ \ \int_0^t ds \left\{  -\sigma^2 \int_{\mathbb{R}} du \, \partial_u G(s,u) \partial_u \tilde{\rho}(s,u) \right. + } \nn \\
& & \left. + \int_{\mathbb{R}} du \, \partial_u G(s,u) \Big[ a\Big( s,u + \int_0^t \frac{C(s)}{2} ds \Big) - \frac{C(s)}{2} \Big] \tilde{\rho}(s,u) \right\} .
\end{eqnarray}
Moreover we have by a simple computation that
$$
a\Big( s,u + \int_0^t \frac{C(s)}{2} ds \Big) - \frac{C(s)}{2} = \frac{\gamma(t,u)}{2}
$$
and thus from (\ref{eq:G1}) and (\ref{eq:F2}) we have that \textbf{(b)} holds in the definition of weak solutions of equation (\ref{eq:eqhydRP2}). $\square$

\medskip
%%%%%%%%%%%%%%%%%%%%%%%%%%%%%%%%%%%%%%%%%%%%%%%%%%%%%%%%%%%%%%%%%%%%%%%%%%%%%
\section{The proof of the Hydrodynamic limit of the EPRS}
%%%%%%%%%%%%%%%%%%%%%%%%%%%%%%%%%%%%%%%%%%%%%%%%%%%%%%%%%%%%%%%%%%%%%%%
\label{sec:proof}
\setcounter{equation}{0}

Denote by $\mathcal{M}=\mathcal{M}(\mathbb{R})$, the space of positive Radon measures on $\mathbb{R}$ endowed with the vague topology. Integration of a function $G$ with respect to a measure $\pi$ in $\mathcal{M}$ will be denoted $\< \pi,G\rangle$. To each configuration $\xi \in \Omega$ and each $N\ge 1$ we associate the empirical measure $\pi^N = \pi^N (\xi)$ in $\mathcal{M}$, where
$$
\pi^N = \frac{1}{N} \sum_{z \in \mathbb{Z}} \xi(z) \delta_{z/N} \, .
$$
Let $D([0,T],\mathcal{M})$ denote the space of right continuous functions with left limits on $\mathcal{M}$ endowed with the Skorohod topology. For each probability measure $\mu$ on $\Omega$, denote by $\mathbb{Q}_{\mu}^N$ the probability measure on $D([0,T],\mathcal{M})$ induced by $\mathbb{P}_{\mu}^N$ and the empirical measure $\pi^N$.

Theorem \ref{theorem:hbep} states that the sequence $\mathbb{Q}_{\mu^N}^N$ converges weakly, as $N \ra \oo$, to the probability measure concentrated on the absolutely continuous trajectory $\pi(t,du)=\zeta(t,u)du$ whose density is the unique weak solution of (\ref{eq:eqhyd1}) (see \cite{kipnislandim}). The proof consists of showing tightness of $\mathbb{Q}_{\mu^N}^N$, that all of its limit points are concentrated on absolutely continuous paths which are weak solutions of (\ref{eq:eqhyd1}) and uniqueness of solutions of this equation.

We have already discussed uniqueness of weak solutions of (\ref{eq:eqhyd1}) in Section \ref{sec:weaksolutions}.

Note that all limit points of the sequence $\mathbb{Q}_{\mu}^N$ are concentrated on absolutely continuous measures since the total mass on compact intervals of the empirical measure $\pi^N$ is bounded by the size of the interval plus $1/N$.

In order to show that all limit points of the sequence $\mathbb{Q}_{\mu^{\scriptscriptstyle{N}}}^N$ are concentrated on weak solutions of (\ref{eq:eqhyd1}) we will need the following result:

\begin{lemma}
\label{lemma:weaksol}
For every smooth function $G:[0,T] \times \mathbb{R} \ra \mathbb{R}$ with compact support and $\delta>0$
\begin{eqnarray}
\lefteqn{ \!\!
\limsup_{N \ra \oo} \mathbb{Q}_{\mu^{\scriptscriptstyle{N}}}^N \left[ \sup_{0\le t\le T} \left| \<\pi_t^N,G\rangle - \<\pi_0^N,G\rangle - \int_0^t \{ \< \pi_s^N, \partial_s G \rangle + \right. \right.} \nn \\ 
& & \qquad \qquad \qquad \qquad \qquad \qquad
\left. \left. + \, \sigma^2 \< \pi_s^N,\Delta G \rangle + \< \pi_s^N, a_s \nabla G \rangle  \} ds  \right| > \delta \right] = 0 \, . \nn
\end{eqnarray}
\end{lemma}

\smallskip

We shall divide the proof of Theorem \ref{theorem:hbep} in four parts: We start proving tightness in Section \ref{sec:tightness}. The Section \ref{sec:weaksol} is devoted to prove of Lemma \ref{lemma:weaksol}. From Lemma \ref{lemma:weaksol}, to conclude the proof that all limit points of the sequence $\mathbb{Q}_{\mu^{\scriptscriptstyle{N}}}^N$ are concentrated on weak solutions of (\ref{eq:eqhyd1}), we have to justify an integration by parts to obtain conditions (a) and (b) in the definition of weak solutions of (\ref{eq:eqhyd1}). This is a consequence of an energy estimate which is the content of Section \ref{sec:energyest}.

\medskip
%%%%%%%%%%%%%%%%%%%%%%%%%%%%%%%%%%%%%%%%%%%%%%%%%%%%%%%%%%%%%%%%%%%%%%%%%%%%%
\subsection{Tightness}
%%%%%%%%%%%%%%%%%%%%%%%%%%%%%%%%%%%%%%%%%%%%%%%%%%%%%%%%%%%%%%%%%%%%%%%
\label{sec:tightness}
\setcounter{equation}{0}

The sequence $\mathbb{Q}_{\mu}^N$ is tight in the space of probability measures on $D([0,T],\mathcal{M})$, if for each smooth function with compact support $G:\mathbb{R} \ra \mathbb{R}$, $\<\pi_t^N,G\rangle$ is tight as a random sequence on $D([0,T],\mathbb{R})$. Now fix such a function, denote by $\mathcal{F}_t = \sigma^2 (\pi_s: s\le t)$, $t \ge 0$, the natural filtration on $D([0,T],\mathcal{M})$, and by $\mathcal{T}_T$ the family of $(\mathcal{F}_t)$-stopping times bounded by $T$. According to Aldous \cite{aldous}, see also Chapter 4 in \cite{kipnislandim}, to prove tightness for $\<\pi_t^N,G\rangle$ we have to verify the following two conditions:
\begin{itemize}
\item[(i)] The finite dimensional distributions of $\<\pi_t^N,G\rangle$ are tight;
\item[(ii)] for every $\epsilon >0$ 
\begin{equation}
\label{eq:aldous}
\lim_{\gamma \ra 0} \limsup_{N \ra \oo} \sup_{\tau \in \mathcal{T}_T} \sup_{\theta \le \gamma} \,
\mathbb{P}_{\mu^{\scriptscriptstyle{N}}}^N \left[ |\<\pi^N_\tau,G\rangle - \<\pi^N_{\tau + \theta},G\rangle | > \epsilon \right] = 0 \, .
\end{equation}
\end{itemize}

Condition (i) is a trivial consequence of the fact that the empirical measure has finite total mass on any compact interval. In order to prove condition (ii), Let us first consider for each smooth function $H:[0,T] \times \mathbb{R} \ra \mathbb{R}¨$ the associated $(\mathcal{F}_t)$-martingale  vanishing at the origin
\begin{equation}
\label{eq:martingal}
M_t^{H,N} = \<\pi_t^N,H\rangle - \<\pi_0^N,H\rangle - \int_0^t (\partial_s + \mathcal{L}_s^N) \<\pi_s^N,H\rangle ds \, .
\end{equation}
In (ii) the function $G$ does not depend on the time, however the martingale $M_t^{H,N}$ is defined for functions varying on time, because we will need it later in the proof of lemma \ref{lemma:weaksol}.  An elementary computation shows that $(\partial_s + \mathcal{L}^N_s) \<\pi^N,H\rangle$ is given explicitly by
\begin{equation}
\label{eq:integrand}
\<\pi^N, \partial_s H\rangle +
\<\pi^N, \Delta_N^p H\rangle + \frac{1}{N} \sum_{z \in \mathbb{Z}} 
\left\{ \frac{1}{N} \sum_{x \le z} b(s,x/N) \right\} \nabla_N H(s,z/N) \xi(z) \, ,
\end{equation}
where $\Delta_N^p$ and $\nabla_N$ denote respectively the discrete Laplacian and gradient:
$$
\Delta_N^p H(s,x/N) = \frac{N^2}{2} \sum_{z\in \mathbb{Z}} p(z) \{H(s,(x+z)/N) + H(s,(x-z)/N) - 2H(s,x/N)\} \, ,
$$
$$
\nabla_NH(s,x/N) = N\{H(s,(x+1)/N) - H(s,x/N)\}. 
$$
The quadratic variation of $M_t^{H,N}$ is denoted by $\<M^{H,N}\rangle_t$ and explicitly given by 
\begin{eqnarray}
\label{eq:quavar}
\lefteqn{ \
\int_0^t ds \Big\{ \frac{1}{N^2} \sum_{z\in \mathbb{Z}_+} p(z) \sum_{|x-y| =z} \Big( \sum_{\omega=x \wedge y}^{x \vee y -1} \nabla_N H(s,\omega/N) \Big)^2 \xi^N_s(x)[1-\xi^N_s(y)] + }\nn \\
& & \!\! + \frac{1}{N^3} \sum_z a(s,z/N) \nabla_N H(s,z/N)^2 \xi^N_s(z) \nn \\
& & \!\! + \frac{1}{N^3} \sum_{z\neq w} a(s,z\wedge w/N) \nabla_N H(s,z/N) \nabla_N H(s,w/N) \xi^N_s(z) \xi^N_s(w) \\
& & \!\! - \frac{1}{N^3} \sum_{z\neq w} \Big\{ \sum_{x=z\wedge w +1}^{z\vee w -1} b(s,x/N) \Big\} H(s,z\vee w/N) \nabla_N H(s,z \wedge w/N) \xi^N_s(z) \xi^N_s(w) \Big\} \nn ,
\end{eqnarray}
where $a^N(s,z) = N^{-1} \sum_{x \le z} b(s,x/N) \, .$
Therefore, for a $G$ not depending on time as before 
$$
\<\pi^N_{\tau + \theta},G\rangle - \<\pi^N_\tau,G\rangle = M_{\tau+\theta}^{G,N} - M_{\tau}^{G,N} + \int_\tau^{\tau + \theta} \mathcal{L}_s^N \<\pi_s^N,G\rangle ds\, .
$$
From the previous expression and Chebyshev inequality, (ii) follows from
\begin{equation}
\label{eq:decmart}
\lim_{\gamma \ra 0} \limsup_{N \ra \oo} \sup_{\tau \in \mathcal{T}_T} \sup_{\theta \le \gamma} \,
\mathbb{E}_{\mu^{\scriptscriptstyle{N}}}^N \left[ \left| M_{\tau+\theta}^{G,N} - M_{\tau}^{G,N} \right| \right] = 0 
\end{equation}
and
\begin{equation}
\label{eq:decint}
\lim_{\gamma \ra 0} \limsup_{N \ra \oo} \sup_{\tau \in \mathcal{T}_T} \sup_{\theta \le \gamma} \,
\mathbb{E}_{\mu^{\scriptscriptstyle{N}}}^N \left[ \left| \int_\tau^{\tau + \theta} \mathcal{L}_s^N \<\pi_s^N,G\rangle ds \right| \right] = 0 \, .
\end{equation} 
Now we show (\ref{eq:decmart}) and (\ref{eq:decint}), completing the proof of tightness and finishing this section.

\medskip

\no \textbf{Proof of (\ref{eq:decmart}):} From the optional stopping theorem and the martingale property
$$
\mathbb{E}_{\mu^{\scriptscriptstyle{N}}}^N \left[ (M_{\tau+\theta}^{G,N} - M_{\tau}^{G,N})^2 \right] = \mathbb{E}_{\mu^{\scriptscriptstyle{N}}}^N \left[ \<M^{G,N}\rangle_{\tau+\theta} - \<M^{G,N}\rangle_\tau \right] \, .
$$
Hence, applying formula (\ref{eq:quavar}), by the Taylor expansion for G, we have that 
\begin{equation}
\label{eq:aldous2}
\mathbb{E}_{\mu^{\scriptscriptstyle{N}}}^N \left[ (M_{\tau+\theta}^{G,N} - M_{\tau}^{G,N})^2 \right] \le \frac{C(G,b)}{N} \, \theta,
\end{equation}
and (\ref{eq:decmart}) holds. Indeed, to deal with the last term in (\ref{eq:quavar}), use the fact that $b$ is Riemann integrable together with the inequality
\begin{eqnarray*}
& & \frac{1}{N^3} \sum_{|z|,|w|\le CN \atop z \neq w} \sum_{x=z\wedge w +1}^{z\vee w -1} b(s,x/N) =
\frac{2}{N^3} \sum_{|z|,|w|\le CN \atop z > w} \sum_{w<x<z} b(s,x/N) \\
& & \le \frac{2}{N^2} \sum_{|w|\le CN} \Big\{ \frac{1}{N} \sum_{x=w+1}^{CN} b(s,x/N) \Big\} \, . \qquad \qquad \qquad \qquad \qquad \qquad \qquad \square 
\end{eqnarray*}

\medskip

\no \textbf{Proof of (\ref{eq:decint}):} 
From formula (\ref{eq:integrand}) and the Taylor expansion for $G$ we obtain that
\begin{equation}
\label{eq:aldous1}
\mathbb{E}_{\mu^{\scriptscriptstyle{N}}}^N \left[ \left| \int_\tau^{\tau + \theta} \! \! \mathcal{L}^N_s \<\pi_s^N,G\rangle ds \right| \right] \le C(G,b) \, \theta. \ \square
\end{equation}

\medskip
%%%%%%%%%%%%%%%%%%%%%%%%%%%%%%%%%%%%%%%%%%%%%%%%%%%%%%%%%%%%%%%%%%%%%%%%%%%%%
\subsection{The proof of Lemma (\ref{lemma:weaksol}) }
%%%%%%%%%%%%%%%%%%%%%%%%%%%%%%%%%%%%%%%%%%%%%%%%%%%%%%%%%%%%%%%%%%%%%%%%%%%%%
\label{sec:weaksol}

Let $G:[0,T] \times \mathbb{R} \ra \mathbb{R}$ be a smooth function with compact support. By Doob inequality, for every $\delta > 0$,
$$
\mathbb{P}_{\mu^{\scriptscriptstyle{N}}}^N \left[ \sup_{0\le t \le T} |M_t^{G,N}| \ge \delta \right] \le 4 \delta^{-2} \mathbb{E}_{\mu^{\scriptscriptstyle{N}}}^N \left[ (M_T^{G,N})^2 \right]
= 4 \delta^{-2} \mathbb{E}_{\mu^{\scriptscriptstyle{N}}}^N \left[  \<M^{G,N}\rangle_T \right], 
$$
which, by the explicit formula for the quadratic variation of $M_t^{G,N}$, see the proof of (\ref{eq:decmart}), is bounded by $C(G,b) \, \theta N^{-1}$. Thus, for every $\delta >0$, 
\begin{equation}
\label{eq:supmart}
\lim_{N \ra \oo} \mathbb{P}_{\mu^{\scriptscriptstyle{N}}}^N \left[ \sup_{0\le t \le T} |M_t^{G,N}| \ge \delta \right] =0.
\end{equation} 
Using (\ref{eq:integrand}) to expand the martingale expression in (\ref{eq:martingal}). Since the Taylor expansion gives us that
$$
\sup_{s \in [0,T]}
\left| N [G(s,x+1/N)-G(s,x/N)] - \nabla G(s,x/N) \right| \le \frac{C(G)}{N},
$$
and $\Delta_N^p H(s,x/N)$ is equal to
$$
\frac{1}{2} \sum_{z\in \mathbb{Z}} p(z) \sum_{y=1}^z \sum_{w=1}^z N \{ \nabla_N H \Big(s,\frac{x-y+w}{N}\Big) - \nabla_N H\Big(s,\frac{x-y+w-1}{N}\Big)\} \, ,
$$
we may replace $\Delta_N^p$ and $\nabla_N$ in (\ref{eq:supmart}) by the usual laplacian and gradient. Moreover the exponencial bound on $b$ given at (\ref{eq:exp}) implies that the Riemann sum $N^{-1} \sum_{x \le [uN]} b(s,x/N)$ converges as $N \ra \oo$ to $\int_0^u b(s,v) dv = a(s,v)$ and we can also replace one by the other in (\ref{eq:supmart}). Thus we obtain that
\begin{eqnarray}
\lefteqn{\!\!\!\!\!\!\!\!\!\!\!\!\!\!
\lim_{N \ra \oo} \mathbb{Q}_{\mu^N}^N \left[ \sup_{0 \le t \le T} \left| \<\pi_t^N,G\rangle - \<\pi_0^N,G\rangle - \int_0^t  ds \{ \<\pi_s^N,\partial_s G\rangle + \right. \right.} \nn \\
& & \qquad \qquad \left. \left. + \sigma^2 \<\pi_s^N,\Delta G\rangle +  \< \pi_s^N, a \nabla G \rangle  \} \right| >\delta \right] =0 \nn
\end{eqnarray}
for all $\delta>0$. $\square$

\medskip
%%%%%%%%%%%%%%%%%%%%%%%%%%%%%%%%%%%%%%%%%%%%%%%%%%%%%%%%%%%%%%%%%%%%%%%%%%%%
\subsection{An energy estimate}
%%%%%%%%%%%%%%%%%%%%%%%%%%%%%%%%%%%%%%%%%%%%%%%%%%%%%%%%%%%%%%%%%%%%%%%%%%%%
\label{sec:energyest}

The next result allows an integration by parts in the expression inside the probability in the statement of Lemma \ref{lemma:weaksol}, proving Theorem \ref{theorem:hbep}.

\medskip
\begin{proposition}%%%%%%%%%%%%%%%%%%%%%%%%%%%%%%%%%%%%%%%%%%%%%%%%%%
\label{theorem:orderofint} 
Every limit point of the sequence $\mathbb{Q}^N$ is concentrated on paths $\zeta(t,u)du$ with the property that $\zeta(t,u)$ is absolutely continuous whose derivative $\partial_u \zeta(t,u)$ is in $L^2([0,T] \times \mathbb{R})$. Moreover
$$
\int_0^T ds \int_{\mathbb{R}} du \, H(s,u) \partial_u \zeta(s,u) = 
 - \, \int_0^T ds  \int_{\mathbb{R}} du \, \partial_u H(s,u) \zeta(s,u) \, ,
$$
for all smooth functions $H:[0,T] \times \mathbb{R} \ra \mathbb{R}$ with compact support.
\end{proposition}%%%%%%%%%%%%%%%%%%%%%%%%%%%%%%%%%%%%%%%%%%%%%%%%%%%%

\medskip
To prove the previous theorem we make use of the following energy estimate:
\medskip
\begin{lemma}%%%%%%%%%%%%%%%%%%%%%%%%%%%%%%%%%%%%%%%%%%%%%%%%%%%
\label{lemma:energyestimate}
There exists $K>0$ such that if $\mathbb{Q}^*$ is a limit point of the sequence $\mathbb{Q}^N$ then
$$
\mathbb{E}_{\mathbb{Q}^*} \left[ \sup_H \left\{ \int_0^T ds \int_{\mathbb{R}} du \, \partial_u H(s,u) \zeta(s,u) - 2 \int_0^T ds \int_{\mathbb{R}} du \, H(s,u)^2 \zeta(s,u) \right\} \right] \le K \, ,
$$
where the supremum is taken over all functions $H$ in $C^{0,1}_0 ([0,T] \times \mathbb{R})$.
\end{lemma}%%%%%%%%%%%%%%%%%%%%%%%%%%%%%%%%%%%%%%%%%%%%%%%%%%%%%

\bigskip

\no \textbf{Proof of Proposition \ref{theorem:orderofint}:} Let $\mathbb{Q}^*$ be a limit point of the sequence $\mathbb{Q}^N$. By Lemma \ref{lemma:energyestimate} for $\mathbb{Q}^*$ almost every path $\zeta(t,u)$ there exists $B=B(\zeta)>0$ such that
\begin{equation}
\label{eq:energyestimate}
\int_0^T ds \int_{\mathbb{R}} du \, \partial_u H(s,u) \zeta(s,u)
-  2 \int_0^T ds \int_{\mathbb{R}} du \, H(s,u)^2 \le B ,
\end{equation}
for every $H \in C^{0,1}_0([0,T],\mathbb{R})$. Note that, since $\zeta<1$ we were able to suppress it in the last integrand. Equation (\ref{eq:energyestimate}) implies that 
$$ 
\lambda(H) := \int_0^T ds \int_{\mathbb{R}} du \, \partial_u H(s,u) \zeta(s,u)
$$
is a bounded linear functional on $C^{0,1}_0([0,T],\mathbb{R})$ for the $L^2$-norm. Extending it to a bounded linear functional on $L^2([0,T],\mathbb{R})$, by Riesz Representation Theorem, there exists a $L^2$ function, denoted by $\partial_u \zeta(s,u)$, such that
$$
\lambda(H)= - \int_0^T ds \int_{\mathbb{R}} du \, H(s,u) \partial_u \zeta(s,u) \, ,
$$
for every smooth function $H:[0,T] \times \mathbb{R} \ra \mathbb{R}$ with compact support. $\square$ 

\bigskip

To show Lemma \ref{lemma:energyestimate}, consider the space $C^{0,1}_0 ([0,T] \times \mathbb{R})$ endowed with the norm 
$$
\|H\|_{0,1} = \sum_{n=0}^{\oo} 2^n \{ \|H \, \mathbf{1}\{(n,n+1)\}\|_\oo + \|\partial_u H \, \mathbf{1}\{(n,n+1)\}\|_\oo \} \,.
$$ 

\medskip

\no \textbf{Proof of Lemma \ref{lemma:energyestimate}:} For every $\epsilon>0$, $\delta>0$, $H:\mathbb{R} \ra \mathbb{R}$ smooth function with compact support and $\xi \in \{0,1\}^{\mathbb{Z}^*}$, denote by $W_N(\epsilon,\delta,H,\xi)$ the following expression
$$
\sum_{x=-\oo}^\oo H(x/N) \frac{1}{\epsilon N} \{\xi^{\delta N}(x) - \xi^{\delta N}(x+[\epsilon N])\} - \frac{2}{N} \sum_{x=-\oo}^\oo H(x/N)^2 \frac{1}{\epsilon N} \sum_{y=0}^{[\epsilon N]} \xi^{\delta N}(x+y) \, ,
$$
where
$
\xi^{\delta}(x) = \delta^{-1} \sum_{y=x}^{x+\delta} \xi(y) .
$
We claim that there exists $K>0$ such that for any dense subset $\{H_l:l\ge 1\}$ of $C^{0,1}_0 ([0,T] \times \mathbb{R})$,
\begin{equation}
\label{eq:A}
\lim_{\delta \ra 0} \lim_{N \ra \oo} \mathbb{E}_{\mu^{\scriptscriptstyle{N}}}^N \left[ \max_{1\le i\le k} \left\{ \int_0^T ds W_N(\epsilon,\delta,H_i(s,\cdot),\xi_s) \right\} \right] \le K \, ,
\end{equation}
for every $k\ge 1$ and every $\epsilon >0$. We postpone the proof of (\ref{eq:A}) and we use it in the sequence. Since $\mathbb{Q}^*$ is a weak limit point of the sequence $\mathbb{Q}_N$, which is concentrated on absolutely continuous trajectories, it follows that
\begin{eqnarray}
\label{eq:A1}
\lefteqn{\limsup_{\delta \ra 0} E_{\mathbb{Q}^*} \left[ \max_{1\le i\le k} \left\{  \int_0^T ds \int_{\mathbb{R}} du \right. \right. } \nn \\
& & \left\{ \epsilon^{-1} H_i(s,u) \left( \delta^{-1} \int_u^{u+\delta} \zeta_s \, dv - \delta^{-1} \int_{u + \epsilon}^{u+\epsilon+\delta} \zeta_s \, dv \right) \right. - \\ 
& & \qquad \quad \left. \left. \left. - 2 \epsilon^{-1} H_i(s,u)^2 \int_u^{u+\epsilon} dv \left( \delta^{-1} \int_v^{v+\delta} \zeta_s \, dv^\prime \right) \right\} \right\} \right] \le K \, , \nn
\end{eqnarray}
for every $k\ge 1$. Since,
$$
\epsilon^{-1} \int_{\mathbb{R}} du H(u) \left\{ \delta^{-1} \int_u^{u+\delta} \zeta_s \, dv - \delta^{-1} \int_{u+\epsilon}^{u+\epsilon+\delta} \zeta_s \, dv \right\} 
$$
is equal to
$$
\int_{-\oo}^\oo du \left\{ \frac{H(u)-H(u-\epsilon)}{\epsilon} \right\} \left\{ \delta^{-1} \int_u^{u+\delta} \zeta_s \, dv \right\} \, ,
$$
letting $\delta \ra 0$ and then $\epsilon \ra 0$, it follows from (\ref{eq:A1}) that
$$
E_{\mathbb{Q}^*} \left[ \max_{1\le i \le k} \int_0^T ds \int_{\mathbb{R}} du \left\{ \partial_uH_i(s,u) \zeta(s,u) -2 H_i(s,u)^2 \zeta(s,u) \right\} \right] \le K \, .
$$
To conclude the proof we apply the monotone convergence theorem for $k\ra \oo$ to replace the maximum over $1\le i\le k$ by the maximum over $i \in \mathbb{N}$ in the integrand above, then note that
$$
\int_0^T ds \int_{\mathbb{R}} du \left\{ \partial_uH(s,u) \zeta(s,u) -2 H(s,u)^2 \zeta(s,u) \right\} ,
$$
is continuous as a real function on $(C^{0,1}_0 ([0,T] \times \mathbb{R}),\| \cdot \|_{0,1})$. $\square$

\bigskip

\no \textbf{Proof of (\ref{eq:A}):} Since $H$ is a continuous function, an integration by parts justify the replacement of $W_N(\epsilon,\delta,H,\xi)$ as $\delta \ra 0$ in (\ref{eq:A}) by
\begin{equation}
\label{eq:W}
\sum_{x=-\oo}^\oo H(s,x/N) \frac{1}{\epsilon N} \{\xi(x) - \xi(x+[\epsilon N])\} - \frac{2}{N} \sum_{x=-\oo}^\oo H(s,x/N)^2 \frac{1}{\epsilon N} \sum_{y=1}^{[\epsilon N]} \xi(x+y) \, ,
\end{equation}
which we denote by $W_N(\epsilon,H,\xi)$. By the entropy inequality 
$$
\mathbb{E}_{\mu^{\scriptscriptstyle{N}}}^N \left[ \max_{1\le i\le k} \left\{ \int_0^T ds W_N(\epsilon,H_i(s,\cdot),\xi_s) \right\} \right]
$$
is bounded by
$$
\frac{H(\mu^N|\nu_{\alpha})}{N} + \frac{1}{N} \log \mathbb{E}^N_{\nu_{\alpha}} \left[ \exp \left\{ N \max_{1\le i\le k} \int_0^T ds W_N(\epsilon,H_i(s,\cdot),\xi_s)  \right\} \right] .
$$
The bound on the entropy given in the statement of Theorem \ref{theorem:hbep} and the elementary inequality $e^{\max a_i} \le \sum e^{a_i}$ imply that this last expression is bounded by
$$
C + \frac{1}{N} \log \mathbb{E}^N_{\nu_{\alpha}} \left[ \sum_{i=1}^k \exp \left\{ N \int_0^T ds W_N(\epsilon,H_i(s,\cdot),\xi_s)  \right\} \right] .
$$ 
Here, since $\max \{ \limsup_{N} N^{-1} \log a_N , \limsup_{N} N^{-1} \log b_N \}$ is greater or equal to $\limsup_{N} N^{-1} \log \{ a_N + b_N \}$, the second term is dominated by
\begin{equation}
\label{eq:W1}
\max_{1\le i\le k} \limsup_{N\ra \oo} \frac{1}{N} \log \mathbb{E}^N_{\nu_{\alpha}} \left[ \exp \left\{ N \int_0^T ds W_N(\epsilon,H_i(s,\cdot),\xi_s)  \right\} \right] .
\end{equation}
We have that the previous expression is bounded by
\begin{equation}
\label{eq:W2}
\max_{1\le i\le k} \int_0^T ds \sup_f \left\{ \int W_N(\epsilon,H_i(s,\cdot),\xi)f(\xi) \nu_{\alpha}(d\xi) - N D(f) \right\} + \frac{1}{N}  \sum_{x \in \mathbb{Z}} \frac{b(s,x/N)}{2 (1-\alpha)}  \, ,
\end{equation}
where the supremum is taken over all cylindrical densities $f$ with respect to $\nu_{\alpha}$, and
\begin{eqnarray*}
D(f)& = &- \int \sqrt{f} L \sqrt{f} d\nu_\alpha \\
& = & \frac{1}{2} \sum_{x,z \in \mathbb{Z}} p(z) \int \xi(x) [1-\xi(x+z)] [\sqrt{f(\xi^{x,x+z}})-\sqrt{f(\xi)}]^2 d\nu_\alpha
\end{eqnarray*}
is the Dirichlet form evaluated on $f$. We assume (\ref{eq:W2}) and show it just after the proof of (\ref{eq:A}). We just have to estimate the first term in (\ref{eq:W2}) whose limsup as $N \ra \oo$ we are now going to show that is in fact non-positive. Write
\begin{eqnarray*}
\lefteqn{\int \{\xi(x) - \xi(x+[\epsilon N])\} f(\xi) \nu_{\alpha}(d\xi) = } \\
& & = 2 \sum_{z \ge 1} \, p(z) \! \sum_{j=0}^{[z^{-1}\epsilon N]-1} \int \{\xi(x+jz) - \xi(x+(j+1)z)\} f(\xi) \nu_{\alpha}(d\xi) \\
& & \quad + \ 2 \sum_{z \ge 1} \, p(z) \! \int \{\xi(x+[z^{-1}\epsilon N]z) - \xi(x+[\epsilon N])\} f(\xi) \nu_{\alpha}(d\xi)
\end{eqnarray*}
where the first term in the right hand side can be rewritten as 
$$
2 \sum_{z \ge 1} \, p(z) \!\!\! \sum_{j=0}^{[z^{-1}\epsilon N]-1} \!\!\! \int \xi(x+(j+1)z) [1-\xi(x+jz)] \{ f(\xi^{x+jz,x+(j+1)z}) - f(\xi)\} \nu_{\alpha} (d\xi).
$$
Hence, we have that
$$
\sum_{x=-\oo}^\oo H(s,x/N) \int  \{\xi(x) - \xi(x+[\epsilon N])\} f(\xi) \nu_{\alpha}(d\xi)
$$
is, for each $B>0$, bounded above by
\begin{eqnarray}
\label{eq:ee1}
\lefteqn{ \sum_{x=-\oo}^\oo B \sum_{z \ge 1} \, p(z) \!\!\! \sum_{j=0}^{[z^{-1}\epsilon N]-1} \Xi(-,x,j,f) } \nn \\
& & \!\!\! + \sum_{x=-\oo}^\oo \frac{H(s,x/N)^2}{B} \sum_{z \ge 1} \, p(z) \!\!\! \sum_{j=0}^{[z^{-1}\epsilon N]-1} \Xi(+,x,j,f) + O(1),
\end{eqnarray}
where $\Xi(\pm,x,j,f)$ denotes the integral
$$
\int \xi(x+(j+1)z) [1-\xi(x+jz)] \left\{ \sqrt{f(\xi^{x+jz,x+(j+1)z})} \pm \sqrt{f(\xi)} \right\}^2 \nu_{\alpha} (d\xi).
$$
Multiplying (\ref{eq:ee1}) by $(\epsilon N)^{-1}$, we bound it above by
$$
B D(f) + \frac{2}{B} \sum_{x=-\oo}^\oo H(s,x/N)^2 \frac{1}{\epsilon N} \sum_{y=0}^{[\epsilon N]} \int \xi(x+y) f(\xi) \nu_{\alpha}(d\xi) + O((\epsilon N)^{-1}) . 
$$
Now taking $B=N$ we obtain from (\ref{eq:W}) that
$$
\int W_N(\epsilon,H(s,\cdot),\xi)f(\xi) \nu_{\alpha}(d\xi) \le N D(f) + O((\epsilon N)^{-1}) \, .
$$
Thus the limsup as $N \ra \oo$ of the first term in (\ref{eq:W2}) is non-positive and (\ref{eq:W1}) is bounded by 
$$\sup_N \frac{1}{N} \sum_{x \in \mathbb{Z}} \frac{b(s,x/N)}{2 (1-\alpha)}$$ 
which is finite since $b$ is integrable. $\square$

\bigskip 
\no \textbf{Proof of (\ref{eq:W2}):} Let $U_{\epsilon}^{N}(s,\xi_s) := W_N (\epsilon,H_i(s,\cdot),\xi_s)$. We wish to estimate
\begin{equation}
\label{eq:entest}
\frac{1}{N} \log \mathbb{E}_{\nu_{\alpha}}^N \left[ \exp \left\{ \int_0^t N U_{\epsilon}^{N}(s,\xi_s) ds \right\} \right] \, .
\end{equation}
Define
$$
(P_{s,t}^N f)(\eta,\xi) = \mathbb{E}_{\xi}^N \left[ 
f(\xi_{t-s}) \exp \left\{ \int_0^{t-s} N U^N_\epsilon(s+r, \xi_r) dr \right\}  \right]
$$
for every bounded real function $f$ on $\{0,1\}^{\mathbb{Z}}$. We have that
\begin{equation}
\label{eq:U1}
\mathbb{E}_{\nu_{\alpha}}^N \left[ \exp \left\{ \int_0^t N U_{\epsilon}^{N}(s,\xi_s) ds \right\} \right] = \int P_{0,t}^N 1 \ d\nu_{\alpha} \le \left\{ \int (P^N_{0,t}1)^2 \ d\nu_{\alpha} \right\}^{\frac{1}{2}} \, .
\end{equation}
We show below that
\begin{eqnarray}
\label{eq:U2}
\partial_s \int (P^N_{s,t}1)^2 \ d\nu_{\alpha} \le
\left\{
2 \sup_f \Big\{ \int N U_{\epsilon}^{N}(s,\xi) f(\xi) \nu_{\alpha}(d\xi) - N^2 D(f) \Big\} \right. + & & \nn \\
+ \left. \sum_{x \in \mathbb{Z}} \frac{b(s,x/N)}{(1-\alpha)} \right\} \int (P_{s,t}^N 1)^2 \, d\nu_{\alpha} & &
\end{eqnarray}
From (\ref{eq:U1}), this inequality and the Gronwall inequality, we have (\ref{eq:W2}).
\smallskip

To complete the proof we show (\ref{eq:U2}). By Feynman-Kac formula, $\{ P_{s,t}: 0\le s\le t\}$ is a semigroup of operators associated to the non-homogeneous generator $N^2 \mathcal{L}^N_s + N U_{\epsilon}^{N} (s,\cdot)$. Moreover, the first Chapman-Kolmogorov equation holds: $\partial_s P_{s,t} = - (N^2 \mathcal{L}^N_s + N U_{\epsilon}^{N} (s,\cdot)) P_{s,t}$. Hence
\begin{eqnarray}
\lefteqn{ - \frac{1}{2} \partial_s  \int (P_{s,t}^N 1)^2 \, d\nu_{\alpha} = \int (N^2 \mathcal{L}^N_s + N U_{\epsilon}^{N} (s,\cdot)) (P_{s,t}^N 1) \, P_{s,t}^N 1 \  d\nu_{\alpha} } \nn \\
& & = \int N^2 L_r^{N,s} (P_{s,t}^N 1) \, P_{s,t}^N 1 \ d\nu_{\alpha} + \int \{ N^2 L + N U_{\epsilon}^{N} (s,\cdot)\} (P_{s,t}^N 1) \, P_{s,t}^N 1 \ d\nu_{\alpha,\gamma}^N \, . \nn
\end{eqnarray}
We shall estimate separately each term in this expression.

\medskip
\no \textbf{Claim 1:} 
$$
\int L_r^{N,s} (P_{s,t}^N 1) \, P_{s,t}^N 1 \ d\nu_{\alpha} \le \sum_{x \in \mathbb{Z}} \frac{b(s,x/N)}{2 (1-\alpha)} \int (P_{s,t}^N 1)^2 \, d\nu_{\alpha} \, .
$$
\no \textbf{Proof of Claim 1:} Denote $P_{s,t}^N 1$ by h. We are considering an integral of the form
$$
\sum_{x \in \mathbb{Z}} \int b(s,x/N) [h(\tau_x \xi) - h(\xi)] \, h(\xi) \ d\nu_{\alpha} \, ,
$$
where $\tau_x$ are defined in Section \ref{sec:serp}. Add and subtract to this the term
$$
\sum_{x\in \mathbb{Z}} \frac{1}{2} \int b(s,x/N) h^2(\tau_x \xi) \ d\nu_{\alpha} = \sum_{x \in \mathbb{Z}} \frac{b(s,x/N)}{2 (1-\alpha)} \int (1-\xi(x)) h^2 \  d\nu_{\alpha} \, .
$$
We end up with two terms, the first is 
$$
- \frac{1}{2} \sum_{x \in \mathbb{Z}}  b(s,x/N) \int [h(\tau_x \xi) - h(\xi)]^2 \ d\nu_{\alpha} \, 
- \frac{1}{2} \sum_{x \in \mathbb{Z}}  b(s,x/N) \int \, h^2 \, d\nu_{\alpha} \, ,
$$
which is negative and may be neglected, and the other is
$$
\sum_{x \in \mathbb{Z}} \frac{b(s,x/N)}{2 (1-\alpha)} \int (1-\xi(x)) h^2 \ d\nu_{\alpha}. \quad \square
$$

\medskip
\no \textbf{Claim 2:} 
\begin{eqnarray}
\lefteqn{\int \{ N^2L + N U_{\epsilon}^{N} (s,\cdot)\} (P_{s,t}^N 1) \, P_{s,t}^N 1 \ d\nu_{\alpha} \le } \nn \\
& & \le \, \sup_f \left\{ \int N U_{\epsilon}^{N}(s,\xi) f(\xi) \nu_{\alpha}(d\xi) - N^2 D(f) \right\}
\int (P_{s,t}^N 1)^2 \, d\nu_{\alpha} \, . \nn
\end{eqnarray}
\no \textbf{Proof of Claim 2:} We have the bound 
$$
\int \{ N^2L + N U_{\epsilon}^{N} (s,\cdot)\} (P_{s,t}^N 1) \, P_{s,t}^N 1 \ d\nu_{\alpha}
\le \Gamma_s^N \int (P_{s,t}^N 1)^2 \ d\nu_{\alpha} \, ,
$$
where $\Gamma_s^N$ is the greatest eigenvalue of the generator $N^2 L + N U_{\epsilon}^{N}(s, \cdot )$. By the variational formula (see Appendix 3 in \cite{kipnislandim}) $\Gamma_s^N$ is equal to
\begin{equation}
\label{eq:varform}
\sup_f \left\{ \int N U_{\epsilon}^{N}(s,\xi) f(\xi) \nu_{\alpha}(d\xi) - N^2 D(f) \right\}. \quad \square
\end{equation}
\smallskip
From the claims 1 and 2 we have (\ref{eq:U2}). $\square$

\bigskip

%%%%%%%%%%%%%%%%%%%%%%%%%%%%%%%%%%%%%%%%%%%%%%%%%%%%%%%%%%%%%%%%%%%%%%%%%%

\end{document}